\documentstyle[12pt]{article}
\newcommand{\const}{\mathop{\rm const}\limits}

\newcommand{\mod}{\mathop{\rm mod}\limits}

\newcommand{\supp}{\mathop{\rm supp}\limits}

\begin{document}

\begin{center}

{\bf CAN THE TAIL FOR MAXIMUM OF CONTINUOUS \\

\vspace{3mm}

RANDOM FIELD BE SIGNIFICANTLY MORE \\

\vspace{3mm}

HEAVY THAN MAXIMUM OF TAILS?  }\\

\vspace{3mm}

\vspace{4mm}

                Eugene Ostrovsky, Leonid Sirota.  \\

\vspace{2mm}
Department of  Mathematic, Bar-Ilan University,  Ramat Gan,
52900, Israel, \\
 e-mails: \ eugostrovsky@list.ru; \ sirota3@bezeqint.net \\

\vspace{4mm}
                    {\sc Abstract.}\\

 \end{center}

 \vspace{3mm}

  \ We construct an example of a continuous centered random process with light tails of finite-dimensional  distribution
but with (relatively) heavy tail of maximum distribution.\par
 \ The apparatus for tails comparison  are embedding results for Orlicz and Grand Lebesgue Spaces (GLS). \par

 \vspace{3mm}


\vspace{3mm}


\vspace{3mm}

\section{Notations. Statement of problem.}

 \vspace{3mm}

   \  The following hypothesis $ H $ has been formulated in an article \cite{Ostrovsky1}, 2008 year:  \\

 \ "Let  $ \theta = \theta(t), \ t \in T  $  be arbitrary separable
random field, centered: $  {\bf E} \theta(t) = 0, $  bounded with probability one:
$ \sup_{t \in T} |\theta(t)| < \infty \ (\mod {\bf P }), $ moreover, may be continuous, if  the set $ T $
is   compact metric space relative some distance. \par
Assume in addition that for some Young (or Young-Orlicz) function $  \Phi(\cdot) $ and correspondent Orlicz norm $ ||\cdot||Or(\Phi) $

$$
\sup_{t \in T}||\theta(t)||Or(\Phi) < \infty. \eqno(1.1)
$$

 Recall that the Luxemburg norm $ ||\xi||Or(\Phi)  $ of a r.v. (measurable function) $  \xi $  is defined as follows:

 $$
  ||\xi||Or(\Phi) = \inf_{k, k > 0} \left\{ \int_{\Omega} \Phi(|\xi(\omega)|/k) \ {\bf P}(d \omega) \le 1 \right\}.
 $$

 \ The Young function $ \Phi(\cdot) $ is by definition arbitrary even convex continuous strictly increasing on the
non-negative right-hand  semi-axis  such that
$$
\Phi(0) = 0, \ \lim_{u \to \infty} \Phi(u) = \infty.
$$

 \ Let also $ \Psi(\cdot) $ be {\it arbitrary} another Young function such that $ \lim_{u \to \infty} \Psi(u) = \infty, \ \Psi <<  \Phi,  $
which denotes by definition

$$
\forall \lambda > 0  \ \Rightarrow \lim_{u \to \infty} \frac{\Psi(\lambda u)}{\Phi(u)} = 0,  \eqno(1.2)
$$
see \cite{Rao1}, p.16. \par

 \ The relation (1.2) is named in the theory of Orlicz spaces as follows: \\
 "the function $  \Psi(\cdot) $ is significantly weaker as  $  \Psi(\cdot), $ write" $ \Psi << \Phi. $  \par

 \ Recall that $ \Psi << \Phi $ implies in particular that the (unit) ball in the space $ Or(\Psi) $ is precompact set
in the space $ Or(\Phi). $ \par

\vspace{4mm}

{\it Open question: there holds (or not)}"

$$
|| \sup_{t \in T} \ |\theta(t)| \ ||Or(\Psi) < \infty. \eqno(1.3)
$$

\vspace{4mm}

 \ The conclusion (1.3) is true for the centered (separable) Gaussian fields \cite{Fernique1},
if the field $ \theta(\cdot) $
satisfies the so-called entropy or generic chaining condition \cite{Ostrovsky2},
\cite{Ostrovsky3},  \cite{Ostrovsky2}, \cite{Talagrand1}, \cite{Talagrand2}, \cite{Talagrand3},
\cite{Talagrand4}; in the case when  $ \theta(\cdot) $ belongs to the domain of attraction of
Law of Iterated Logarithm \cite{Ostrovsky5} etc.\par

 \ Notice that if the field $ \theta(t) $ is continuous $ ( \mod {\bf P}) $  and satisfies the condition
(1.1), then {\it there exists } an Young function $ \Psi(\cdot), \ \Psi(\cdot) << \Phi(\cdot) $
for which the inequality  holds, see \cite{Ostrovsky1}. \par

 \ The condition of a form $ ||\xi||Or(\Psi) < \infty $ described the tail behavior for
the distribution of the random  variable $ \xi. $  Another approach which was used in the
monograph M.Ledoux and M.Talagrand \cite{Talagrand1}, p. 309-317 is related in fact with generalized
Lorentz  (more exactly, Lorentz-Zygmund) norm $ ||  \ \xi \ ||L(v):  $

$$
 || \xi ||L(v) \stackrel{def}{=} \sup_{ A: {\bf P}(A)> 0}  \left[ \frac{1}{v({\bf P}(A))} \cdot
 \int_A |\xi(\omega)| \ {\bf P}( d \omega) \right].
$$
 \ Here $ v = v(z), \ z \in (0,1] $ is continuous strictly monotonically increasing function   such that
 $ v(0) = v(0+) = 0. $ \par

 \ Notice that in all this cases  the inequality (1.3) is true with replacing the function $  \Psi $
on the function $ \Phi.$\par

\vspace{3mm}

 \ The negative answer on the problem  (1.3) for the Orlicz spaces  was given in the preprint \cite{Ostrovsky10},
 especially in the case of infinite "probability" measure  $  {\bf P}. $ \par

 \vspace{3mm}

  \ {\bf  Our target in this short report is to extend a negative answer on the formulated above  hypothesis also on the case the
 so-called Grand Lebesgue Spaces (GLS),   by means of construction of correspondent counterexamples. }\\

\vspace{3mm}

 \ Note that the majority of Orlicz's spaces are a particular or extremal cases of the Grand Lebesgue Spaces,
for example, exponential Orlicz's spaces, classical Lebesgue-Riesz spaces $  L(p) $ etc., see
  \cite{Ostrovsky3},  \cite{Ostrovsky8}.\par

\vspace{4mm}

\section{ Several notations and definitions. Auxiliary facts.}\par

\vspace{3mm}

 {\bf A.} A triplet $ (\Omega, \cal{B}, {\bf P} ),  $  where $ \Omega = \{\omega\} = \{x\}  $ is arbitrary set,
$ \cal{B} $ is non-trivial sigma-algebra  subsets $ \Omega  $ and  $ {\bf P}  $ is non-zero non-negative completely
additive measure defined on the $  \cal{B} $ is called a probabilistic space, even in the case when $ {\bf P}(\Omega) = \infty. $ \par
 \  We denote as usually for the random variable $ \xi $ (r.v.) (i.e. measurable function $ \ \xi: \Omega \to R )  $

 $$
 |\xi|_p = [ {\bf E} |\xi|^p ]^{1/p} = \left[ \int_{ \Omega  } |\xi(\omega)|^p \ {\bf P}(d \omega)  \right]^{1/p}, \ p \ge 1;
 $$

$$
L_p = \{ \xi, \ |\xi|_p < \infty  \}. \eqno(2.0)
$$

\vspace{3mm}

{\bf B.} The so-called Grand Lebesgue Space $ G \psi = G \psi(a,b) $  with norm $  ||\cdot ||G\psi   $ is defined (in this article) as follows:

$$
G \psi = \{\xi, \ ||\xi||G \psi < \infty \}, \  ||\xi||G \psi \stackrel{def}{=} \sup_{ p \ge 1 } \left[ \frac{|\xi|_p}{\psi(p)} \right].
\eqno(2.1)
$$
  \ Here $ \psi = \psi(p) $ is some continuous  function  defined on some numerical interval $ p \in (a,b), $ where
$  1 \le a < b \le \infty,  $  and  such that $ \inf_{p \in (a,b)} \psi(p) > 0. $\par

 \ We will denote

 $$
 (a,b) = \supp \psi(\cdot)
 $$
and define formally $ \psi(p) = + \infty, \ p \notin [a,b]. $ \par

 \ The set of all such a functions will be denotes by $  \Psi(a,b);  $ define also

$$
\Psi \stackrel{def}{=} \cup_{1 \le a < b \le \infty} \Psi(a,b).
$$

 \ The detail investigation of this spaces (and more general spaces) see in \cite{Liflyand1},  \cite{Ostrovsky8}. See also
 \cite{Fiorenza1}, \cite{Fiorenza2}, \cite{Iwaniec1}, \cite{Iwaniec2}, \cite{Kozachenko1} etc. \par

\vspace{4mm}

 \ {\bf Example 2.1.}\par

 \   An important for us fact about considered here spaces is proved in \cite{Ostrovsky6}: if $  {\bf P}(\Omega) = 1 $ and
 $ a = 1, \ b = \infty, $ then under some simple conditions   they coincide with certain exponential  Orlicz's spaces
  $ Or(\Phi_{\psi}). $ For instance, if $ {\bf P}(\Omega) = 1 $ and  $ \psi(p) = \psi_{1/2}(p): =\sqrt{p}, $ then the space
  $ G\psi_{1/2} $ consists on all the subgaussian (non-centered, in general case) r.v. $ Or(\Phi_{\psi_{1/2}})  $
   for which $ \Phi_{\psi_{1/2}}(u) = \exp(u^2/2) - 1.  $  \par
 The Gaussian distributed r.v. $ \eta $ belongs to this space.  Another example: let $ \Omega = (0,1) $ with usually Lebesgue measure and

 $$
 f_{1/2}(\omega) = \sqrt{|\log \omega}|, \ \omega > 0; \ f_{1/2}(0)= 0.
 $$
  It is easy to calculate using Stirling's formula for the Gamma function:
 $$
 | f_{1/2}|_p \asymp \sqrt{p}, \ p \in (1,\infty).
 $$
  The tail behavior:

  $$
  {\bf P} ( f_{1/2} > u ) = \exp (-u^2).
  $$

\vspace{3mm}

{\bf Example 2.2.} \par

 \ If we define the {\it degenerate } $ \psi_{(r)}(p), \ r = \const \ge 1 $ function as follows:
$$
\psi_{(r)}(p) = \infty, \ p \ne r; \hspace{3mm} \psi_{(r)}(r) = 1
$$
and agree $ C/\infty = 0, C = \const > 0, $ then the $ G\psi_r(\cdot) $ space coincides
with the classical Lebesgue space $ L_r. $ \par

\vspace{3mm}

 \ {\bf Example 2.3.} \ An used further example:

 $$
 \psi^{(\beta,b)}(p) = (b-p)^{-\beta}, \ 1 \le p < b, \ \beta = \const \ge 0; \ G_{\beta,b}(p) := G_b\psi^{(\beta,b)}(p). \eqno(2.2)
 $$

\vspace{4mm}

{\bf C.} Recall that sets $ A_1, A_2, \ A_i \in \cal{B} $ are disjoint,  if $ A_1\cap A_2 = \emptyset. $ The sequence of a functions
$  \{h_n \}, n =1,2,3  \ldots  $ is said to be {\it disjoint}, or more exactly {\it pairwise disjoint, if }

$$
\forall i,j; i \ne j \ \Rightarrow \ h_i \cdot h_j = 0. \eqno(2.3)
$$
 \ If the sequence of a functions $  \{h_n \} $  is pairwise disjoint, then

$$
|\sum_n h_n|_p^p = \sum_n |h_n|_p^p, \hspace{6mm}  \sup_n |h_n(x)| = \sum_n |h_n(x)|. \eqno(2.4)
$$

\vspace{4mm}

{\bf D.} \ We denote as ordinary  for any measurable set $ A, \ A \in \cal{B} $ it indicator function by
$ I(A) = I_A(\omega).   $\par

\vspace{4mm}

{\bf E.} \ Let $  \phi = \phi(p) $ and $  \psi = \psi(p), \ p \in (a,b), \ 1 \le a < b \le \infty $  be two
functions  from one and the same $  G\psi $ space $  \Psi(a,b).  $ By definition, see \cite{Ostrovsky8}, \cite{Liflyand1},
the function $  \phi(\cdot) $ is significantly weaker in the sense of Grand Lebesgue Spaces,
as one $  \psi(\cdot), $ write also $ \phi << \psi, $  iff (attention, please!)

$$
\lim_{\phi(p) \to \infty} \frac{\phi(p)}{\psi(p)} = \infty. \eqno(2.5)
$$

 \ The relation (2.5) is simpler in comparison with ones in (1.2).

 \ As before, the relation $ \phi << \psi $ in the GLS sense
 implies in particular that the (unit) ball in the space $ G\psi $ is precompact set
in the space $ G\phi; $ on the other words,  compact embedding. \par

\vspace{4mm}

\section{Main result.}

\vspace{4mm}

 \ {\bf Theorem 3.1.} {\it The proposition of hypothesis $ H $  is not true even in the Grand Lebesgue Spaces. } \par

 \ In detail, there exist: \par

\vspace{3mm}

 \ {\bf A.}  A non-trivial $ \psi(\cdot) $ function  from the set $  \Psi(1,b), \ 1 < b \le \infty $
and  {\it compact } non-trivial metric space $ (T,d) = (\{t \}, d). $\par

\vspace{3mm}

 \ {\bf B.} A centered and continuous in the $ G\psi $ sense

$$
\forall s \in T \ \Rightarrow \lim_{t \to s} ||\theta(t) - \theta(s)||G\psi = 0 \eqno(3.1)
$$
and with probability one

$$
{\bf P} (  \theta(\cdot) \in C(T,d) ) = 1 \eqno(3.2)
$$
numerical valued random process (field) $ \theta = \theta(t) = \theta(t,\omega) $
defined aside from the probabilistic space on  our metric space $ (T,d) = (\{t \}, d), $
such that

$$
\sup_{t \in T} ||\theta(t)||G\psi < \infty. \eqno(3.3)
$$

\vspace{3mm}

{\bf C.} A $  \Psi(1,b)  $ function $ \phi = \phi(p) $  which is significantly weaker in the Grand
Lebesgue Space sense as the function $  \psi: \  \phi << \psi $   but herewith

\vspace{3mm}

{\bf D.}

$$
  || \sup_{t \in T} \ |\theta(t)| \ ||G\phi = \infty. \eqno(3.4)
$$

{\bf Proof.}\\

\vspace{3mm}

{\bf 1.} We choose in the sequel  as  the metric space $  (T,d) $ the set of positive integer
numbers with infinite associated point  which we denote by $  \infty: $

$$
T = \{ 1,2,3, \ldots, \infty   \}.  \eqno(3.5)
$$
 The distance $  d $ is defined as follows:
 $$
 d(i,j) = \left| \frac{1}{i} - \frac{1}{j} \right|, \ i,j < \infty; \ d(i,\infty) = d(\infty,i)= \frac{1}{i}, \ i < \infty;
 \eqno(3.5a)
 $$
and obviously $ d(\infty,\infty) = 0.  $\par
 The pair $ (T,d) $ is compact (closed) metric space  and the set $ T $ has an unique limit point $ t_{0} = \infty. $
For instance, $ \lim_{n \to \infty} d(n,\infty) = 0.   $ \par

\vspace{4mm}

{\bf 2.}  It is enough to consider the case $ {\bf P} (\Omega) = 1.  $  More detail, let $ \Omega = (0,1)  $
with ordinary Lebesgue measure. \par

 \ Let    $  f = f(x), \ x \in \Omega = (0,1)  $ be non-zero non-negative integrable function belonging to the
space $  L_4.  $  Define a following $  \psi \ -  $ function:

$$
\nu(p) = |f|_p = \left[  \int_0^1 |f(x)|^p \ dx  \right]^{1/p}, \ 1 \le p \le 4.  \eqno(3.6)
$$
 On the other words, $  \nu(\cdot) $ is a natural function for the function $  f. $ Evidently, $  \nu(\cdot) \in G\psi(1,4).  $\par

  Introduce also  the following numerical sequences

$$
c_n := n^{\beta},  \ \beta = \const > 0, \ n = 2,3, \ldots; \eqno(3.7)
$$

$$
\Delta_n := C(\beta) \cdot n^{ - 4 \beta - 1 }, \ C(\beta): \ \sum_{n=1}^{\infty} \Delta_n = 1; \
a_n = a(n):= \sum_{m=n}^{\infty} \Delta_n;  \eqno(3.8)
$$
and define also sequence of a functions and likewise the following  positive random process
$ \theta(t) = g_n, \ n = t, \ t,n \in T, \ \Omega = \{x \},  $

$$
g_n(x) = c(n) \ f \left( \frac{x-a(n)}{\Delta(n)} \right) \ I_{(a(n+1), a(n)) }(x), \ x \in \Omega, \
g_{\infty}(x) = 0;   \eqno(3.9)
$$

$$
g(x) = \sum_{n=1}^{\infty} g_n(x) =  \sum_{n=1}^{\infty} c_n  \ f \left( \frac{x-a(n)}{\Delta(n)} \right) \ I_{(a(n+1), a(n)) }(x).  \eqno(3.10)
$$

 \ Note that the sequence of r.v.  $ \{ g_n(x) \}  $ consists on  non-negative and  disjoint functions, therefore

 $$
 \sup_n g_n(x) = \sum_n g_n(x) = g(x), \hspace{5mm} |\sup_n g_n|^p_p = \sum_n |g_n|^p_p. \eqno(3.11)
 $$

\vspace{3mm}

  \ Note also that the functions  $ g_n $ are disjoint and following $ \sup_n |g_n(x)| < \infty  $ almost surely. \par
We calculate using the relations (3.7) - (3.11):

$$
|g_n|_p^p = c^p(n) \ \Delta_n \ \nu^p(p) = C(\beta) \ n^{p \beta - 4 \beta - 1} \ \nu^p(p),
\ 1 \le p \le 4,  \eqno(3.12)
$$
therefore

$$
\sup_{p \in [1,4]}  \sup_n |g_n|_p^p \le C(\beta) \ \nu^4(4)  < \infty \eqno(3.13a)
$$
or equivalently

$$
\sup_n |g_n(\cdot)| \in L_4. \eqno(3.13b)
$$

 \ Moreover, $ g_n \to 0  $ almost everywhere.  Indeed, let $ \epsilon $ be arbitrary positive number. We get
applying the estimate (3.12) at the value $ p = 1 $  and Tchebychev's inequality

$$
\sum_n {\bf P} ( |g_n| > \epsilon ) \le C(\beta) \sum_n \frac{n^{ - 3 \beta - 1  }}{\epsilon} < \infty.
$$
 Our conclusion follows immediately from the lemma of Borel-Cantelli. \par

\vspace{3mm}

 \  So, the random process $ \theta(t) = g_n, $ where $ \ n = t \ $ satisfies the condition (1.1)
relative the $ \Psi \ -   $ function $ \psi_{(4)}(p)  $ and is continuous almost everywhere relative the
distance function $  d = d(t,s). $ \par

 \vspace{3mm}

  \  Let now find  the exact up to multiplicative constant expression for the natural function of the r.v.
 $  \sup_n |g_n(x)| $  as $  p \to 4-0.  $  We have:

 $$
  | \ \sup_n |g_n| \ |_p^p =  \sum_n |g_n|_p^p    = \sum_n c^p(n) \ \Delta_n \  \nu^p(p) =
 $$

$$
  = C(\beta) \ \ \nu^p(p)  \ \sum_n n^{p \beta - 4 \beta - 1} \ \sim \frac{C_1(\beta)}{4 - p};
  \eqno(3.14)
$$

$$
| \ \sup_n \ | g_n|  \ |_p  \sim C_2(\beta)(4 - p)^{-1/4}. \eqno(3.15)
$$

\vspace{3mm}

 \  Thus, we can choose for the proposition of theorem 3.1
 as  the $ \Psi \ - $  function $ \psi(p) $  the function  $  \psi_{(4)}(p), $
which is in turn equivalent to the following $   \Psi \ -  $ function

 $$
 \psi^{(0,4)}(p) = 1, \ 1 \le p < 4,
 $$
and correspondingly  to take

$$
\phi_0(p):=  (4 - p)^{-1/8} = \psi^{(1/8, 4)}(p), \ 1 \le p < 4, \eqno(3.16)
$$
see example 2.3. \par
 \ Obviously,

 $$
  \phi_0(\cdot) << \psi^{(0,4)}(\cdot) \eqno(3.17)
 $$
and

$$
|| \ \sup_n \ | g_n|  \ ||G\phi_0 = \infty, \eqno(3.18)
$$

\vspace{3mm}

  \ In order to obtain the centered needed process $ \theta(t)  $  with at the same properties,
 we consider the sequence  $  \tilde{g}_n(x) = \epsilon(n) \cdot g_n(x),  $  where $ \{\epsilon(n) \} $ is a Rademacher
 sequence independent on the $  \{g_n\}: $

$$
{\bf P}(\epsilon(n) = 1) = {\bf P}(\epsilon(n) = -1) = 1/2; \eqno(3.19)
$$
then

$$
 | \ \tilde{g}_n(x) \ | = | \ g_n(x) \ |,  \ | \ \tilde{g}_n \ |_p = | \ g_n \ |_p \eqno(3.20)
$$
 and the sequence $ \{\tilde{g}_n \} $ is also pairwise disjoint (Rademacher's symmetrization).\par
  This completes the proof of our theorem 3.1.

\vspace{4mm}

{\bf Remark 3.1.} The constructed process $ \theta(t) $ give us a new example of centered continuous random process with relatively
 light tails of finite-dimensional distribution, but for which the so-called entropy and generic chains series divergent.\par

\vspace{3mm}

{\bf Remark 3.2.} The proposition of our theorem 3.1 remains true if we use instead the space of continuous function $ C(T,d)  $
arbitrary  separable Banach space.\par

\vspace{3mm}

{\bf Remark 3.3.} Our constructions are likewise to ones in the author's preprint
\cite{Ostrovsky10}.\par

\vspace{4mm}

\section{Concluding remarks.}

\vspace{4mm}

 {\bf General boundedness condition.}\par

\vspace{3mm}

 \ There are many works devoted to deducing of sufficient condition  (entropy conditions as well as
conditions based on the so-called majorizing measure conditions)  for boundedness (continuity)
of the random fields, see e.g. \cite{Fernique1}, \cite{Kozachenko1}, \cite{Ostrovsky3}, \cite{Pizier1},
\cite{Talagrand2}  - \cite{Talagrand5}. \par
 \ Note in addition that if $  \kappa(t), \ t \in T,  $ where $ T  $ is arbitrary set, is separable
numerical random process (field),   and

$$
T = \cup_{k=1}^{\infty} T_k
$$
is countable non-random partition of the set $  T,  $ then

$$
{\bf P} (\sup_{t \it T} \kappa(t) > u) \le \sum_{k=1}^{\infty} {\bf P} (\sup_{t \in T_k} \kappa(t) > u),
\ u = \const, \eqno(4.1)
$$
where each of summands in (4.1)  may be estimated by means of entropy or majorizing measure methods. \par

\end{document}